\documentstyle[12pt]{article}

\newtheorem{prop}{Proposition}
\newtheorem{lem}[prop]{Lemma}
\newtheorem{cor}[prop]{Corollary}

\newcommand{\norm}[1]{\mbox{$\|#1\|$}}

\newcommand{\om}{\omega}
\font\tenBbb=msbm10  \font\sevenBbb=msbm7  \font\fiveBbb=msbm5
\newfam\Bbbfam
\textfont\Bbbfam=\tenBbb \scriptfont\Bbbfam=\sevenBbb
\scriptscriptfont\Bbbfam=\fiveBbb
\def\Bbb{\fam\Bbbfam\tenBbb}
\def\R{{\Bbb R}}

\def\N{{\Bbb N}}

\title{On a question of Haskell P. Rosenthal}
\author{Valentin Ferenczi, Anna Maria Pelczar and
Christian Rosendal}

\begin{document}

\newcommand{\eps}{\epsilon}

\newcommand{\compl}{\complement}
\newcommand {\co} {\complement}
\newcommand {\iso} {\cong}
\newcommand{\isom}{\simeq}

\newcommand {\tom} {\emptyset}
\newcommand{\emb}{\sqsubseteq}
\newcommand{\inj}{\hookrightarrow}
\newcommand{\begr}{\upharpoonright}
\newcommand{\surj}{\twoheadrightarrow}
\newcommand{\bij}{\longleftrightarrow}

\newcommand{\hviss}{\longleftrightarrow}
\newcommand{\hvis}{\Longleftarrow}
\newcommand{\saa}{\Longrightarrow}
\newcommand{\til}{\longrightarrow}

\newcommand{\Lim}[1]{\mathop{\longrightarrow}\limits_{#1}}

\newcommand {\Del}{ \; \Big| \;}
\newcommand {\del}{ \; \big| \;}

\newcommand {\Mgdv}{\Big\{}
\newcommand {\mgdv}{\big\{}

\newcommand {\Mgdh}{\Big\}}
\newcommand {\mgdh}{\big\}}

\newcommand {\Intv}{\Big[}
\newcommand {\intv}{\big[}

\newcommand {\Inth}{\Big]}
\newcommand {\inth}{\big]}

\newcommand {\For}{\Bigcup}
\newcommand {\for}{\bigcup}

\newcommand {\Snit}{\Bigcap}
\newcommand {\snit}{\bigcap}

\newcommand {\og}{\; \land \;}
\newcommand {\eller}{\; \vee\;}
\newcommand{\ikke}{\lnot}

\newcommand {\go} {\mathfrak}
\newcommand {\ku} {\mathcal}
\newcommand {\un} {\underline}

\newcommand {\e} {\exists}
\renewcommand {\a} {\forall}
\newcommand{\fed}{\boldsymbol}

\newcommand{\pf}{

\smallskip

\noindent {\it Proof : }}

\newcommand{\pff}{$\hfill  \Box$

\smallskip }

\newcommand{\pffclaim}{$\hfill \textrm{proof of claim} \Box$

\smallskip }

\maketitle

\begin{abstract}
We consider a normalized basis in a Banach space with the following property: 
any normalized block
sequence of the basis has a subsequence equivalent to the basis. We show that
under uniformity or other natural assumptions, a basis with this property 
is equivalent to the unit 
vector basis of $c_0$ or $\ell_p$. We also address an analogous problem 
concerning spreading models. 
\end{abstract}

Haskell P. Rosenthal has posed the following problem on basic sequences
 in a Banach space:

\paragraph{Problem} {\em Let $X$ have a normalized basis $\{e_i\}$ with the property
 that every normalized block basis admits a subsequence equivalent to
$\{e_i\}$. Is $\{e_i\}$ equivalent to the unit vector
basis of $l_p$ or $c_0$?}

\

Let us recall a well-known theorem of Zippin, which states that a
normalized basis of a Banach space such that all normalized 
block bases are equivalent
(to the original basis) must be equivalent to the unit vector
basis of $l_p$ or $c_0$.

The problem of Rosenthal is of particular
interest as of a "mixed" Ramsey type, in the sense that
it links two types of "subbases" of a given basis: namely subsequences and
block sequences. An instance of a theorem which mixes a property
concerning subsequences and a property concerning block bases was
given by the second named author
in \cite{pel}. She proved that a Banach space saturated with subsymmetric 
sequences must contain a minimal subspace.
 
  Let us notice that this mixing is necessary to make Rosenthal's problem
significant. Indeed, a weakening of the Rosenthal property would be
to assume that every subsequence has a further subsequence
equivalent to the basis. An  application of Ramsey theorem would
give us that the basis is subsymmetric. But
obviously not every subsymmetric basis is 
equivalent to $c_0$ or
$l_p$ (take the basis of Schlumprecht's space for example).

On the other hand we may weaken the Rosenthal property by only
requiring that any block sequence has a further block sequence
equivalent to $\{e_i\}$. Let us call a basis with this property
{\em block equivalence minimal}. The correct setting for such a
property is Gowers' block version of Ramsey theorem. A standard
diagonalization yields that for some constant $C>0$, any block
sequence has a further block sequence which is $C$-equivalent to
$\{e_i\}$, and by Gowers' dichotomy theorem we get the existence
of a winning strategy for Player 2 in Gowers' game to produce
sequences $C+\epsilon$-equivalent to $\{e_i\}$. But again,
Schlumprecht's space is a non-trivial example (see \cite{sch} or
\cite{AS}). Actually, by the proof of Theorem 3.4 in \cite{pel},
 any Banach space
saturated with subsymmetric sequences contains a block equivalence 
minimal basic sequence.

Rosenthal's problem is closely related to a problem of Argyros
concerning spreading models: if all spreading models in a Banach
space are equivalent, must they be equivalent to $c_0$ or $l_p$?
We inspire ourselves from results of Androulakis, Odell,
Schlumprecht and Tomczak-Jaegermann   (\cite{and}) about Argyros'
question, to prove that the answer to Rosenthal's question is
positive if uniformity is assumed (Proposition \ref{prop1}), or when $1$ is
in the Krivine's set of the basis (Corollary \ref{cor6}). We show that the
answer is also positive if $X$ and $X^*$ satisfy the property of
Rosenthal (Proposition \ref{prop7}), or when the selection of the
subsequence in the definition of Rosenthal's property can be
chosen to be continuous (Proposition \ref{prop9}). Let us notice that this
is in contrast with the case of block equivalence minimality,
where uniformity comes as a consequence of the definition, as well
as continuity (because of the previous remark using Gowers'
dichotomy theorem). Finally, after relating Rosenthal's
to Argyros' question, we show how results of descriptive set
 theory may be used to get a dichotomy concerning the number
  of non-equivalent spreading models in a Banach space 
(Proposition \ref{prop10}), and we show that $c_0$ or $l_p$ embeds in $X$ if
there exists a continuous way of picking subsequences generating
spreading models (Proposition \ref{prop11}).

\

Let us give a definition for the main property in this paper. A normalized 
basis $\{e_n\}$
such that any normalized block basis of $\{e_n\}$ has a subsequence which is equivalent
to $\{e_n\}$ will be said to have {\em Rosenthal's} property, or
in short to be a
{\em Rosenthal} basis.

\

First some easy remarks. We notice that a Rosenthal basis must
be subsymmetric. Indeed let ${\ku A}\subset [\om]^{\aleph_0}$ be the set
of subsequences of $\om$
 giving a subsequence of $\{e_i\}$ equivalent to $\{e_i\}$. Then $\ku
A$ is clearly Borel and therefore by the Galvin-Prikry theorem, there is an infinite subset $H$ of $\om$
such that either
 $[H]^{\aleph_0}\subset \ku A$ or  $[H]^{\aleph_0}\subset  {\ku A}^C$.
Evidently the last
 possibility contradicts Rosenthal's property.
So $\{e_i\}_\om \sim \{e_i\}_H$ and $\{e_i\}_H$ is subsymmetric,
hence $\{e_i\}_\om$ also.

By renorming we can assume that $\{e_i\}$ is invariant under
spreading ($1$-equivalent to its subsequences). Now Brunel and
Sucheston have observed that for a normalized basic sequence
$\{t_i\}$, invariant under
spreading, the difference sequence $\{t_{2i+1}-t_{2i}\}$ is
suppression unconditional (ie. the norm decreases as the support
diminishes). By Rosenthal's property ${\{e_{2i+1}-e_{2i}\}}$  is
equivalent to $\{e_i\}$, and is also invariant under spreading. So
we may always assume that a Rosenthal basis is both invariant
under spreading and suppression unconditional. Let us notice at
this point that
 according to
Schlumprecht's terminology, spaces with a Rosenthal basis
 are exactly spaces of Class 1
with a subsymmetric basis. So for once, our favorite non-trivial
example $S$ will not do: it is of Class 2! (see \cite{Sc2}).

\

We fix some notational matters:
we say that a block basis $\{x_i\}$ over $\{e_i\}$ is identically
 distributed if there are scalars $r_0,\ldots,r_k$ and natural numbers
$m_0<m_0+k<m_1<m_1+k<m_2<m_2+k<\ldots$ such that
$x_i=r_0e_{m_i}+\ldots+r_ke_{m_i+k}$.

\

We show that with a uniformity condition added in the hypothesis the answer
is positive. In fact we get a bit more:

\begin{prop}\label{prop1} Let $\{e_i\}$ be a normalized basic sequence
and $K\geq 1$ be
a constant such that any identically distributed normalized  block basis
admits a subsequence $K$-equivalent to $\{e_i\}$. Then
$\{e_i\}$ is equivalent to the unit vector basis of $c_0$ or $l_p$.
\end{prop}

\pf
The proof in the previous remark still goes to show that without
loss of generality, we may assume that
$\{e_i\}$ is both invariant under spreading and suppression
unconditional.
Now under these conditions, Krivine's theorem takes a particularly simple
form:

(Krivine) Let $\{t_i\}$ be a suppression unconditional basis, invariant
under spreading. Then there is
a $p\in [1,\infty]$ such that for all $k<\om,0<\eps$, there are identically
distributed blocks $x_1<x_2<\ldots<x_k$ that are
$(1+\eps)$-equivalent to the unit vector basis of $l_p^k$.

The set of $p$'s satisfying this assertion is called the {\em Krivine
set}. Take a $p$ in this set for our basis $\{e_i\}$, then for any $k$
there is a norm one block $x(k)$, such that
 taking $k$ successive copies of this vector $x_1(k)<x_2(k)<\ldots<x_k(k)$
gives a sequence $2$-equivalent to the unit vector
basis in $l_p^k$.

Taking now infinitely many copies of  $x(k)$:
$x_0(k)<x_1(k)<x_2(k)<\ldots$, it can be observed that
the sequence is identically distributed, so as before it must be
$K$-equivalent to $\{e_i\}$.
But this means that $\{e_0,\ldots,e_{k-1}\}$ is $2K$-equivalent to $l_p^k$,
and as $k$ was arbitrary, $\{e_i\}$ must be $2K$-equivalent to $l_p$ or
$c_0$ if $p=0$.
\pff

\paragraph{Remark} The uniformity condition is necessary in this result.
Indeed take any invariant under spreading $1$-unconditional basis $\{e_i\}$
(like our usual example of
the unit basis of Schlumprecht's space which is not equivalent to $c_0$ or
$l_p$): then any
identically distributed sequence
$\{x_i\}=\{r_0e_{m_i}+\ldots+r_ke_{m_i+k}\}$ is equivalent to $\{e_i\}$, as
proved by
the relation
$$(\max_{j}{|r_j|})\norm{\sum \lambda_i e_{i}} \leq
 \norm{\sum
\lambda_i x_i} \leq
(\sum_{j=0}^k{|r_j|})\norm{\sum \lambda_i e_i},$$
for all  sequences $(\lambda_i)$ in $c_{00}$.

Notice also that this is an opposition to the property of block 
equivalence minimality,
where uniformity is a direct consequence of the property.

\

We now study Rosenthal's problem without the uniformity condition.  First
we notice that the only
relevant case is the reflexive one, as showed by the next lemma. We
need some notation about spreading models. A sequence
$\{x_i\}$ in a Banach space $  X$ is called seminormalized if there are
real numbers
$0<c<C$ such that
$c<\|x_i\|<C,
\; \a i$. Let  $\{x_i\}$ in a Banach space $  X$ be a
seminormalized basic sequence. Suppose that
$$\a r_0,\ldots,r_k\in \R\; \e t\in \R \;\a \eps >0 \;\e N\;\a
N<l_0<\ldots<l_k,$$ $$
\Del\|r_0x_{l_0}+\ldots+r_kx_{l_k}\|-t\Del<\eps$$ (or more
intuitively
$\lim_{l_0<\ldots<l_k,\;l_0\to\infty}\|r_0x_{l_0}+\ldots+r_kx_{l_k}\|$
exists), then we say that $\{x_i\}$ generates a spreading model
$\{\tilde{x_i}\}$ with the norm defined as follows:
$$\|r_0\tilde{x_{0}}+\ldots+r_k\tilde{x_{k}}\|\
:=\lim_{l_0<\ldots<l_k,\;l_0
\rightarrow\infty}\|r_0x_{l_0}+\ldots+r_kx_{l_k}\|$$
The spreading model  $\{\tilde{x_i}\}$ is then a basic
sequence, invariant under spreading.
 Furthermore it is easily seen that the basic constant of $\{\tilde x_i\}$
is majorized by that of $\{x_i\}$.
Moreover any subsequence of  $\{x_i\}$ generates the same spreading model.

\begin{lem}\label{lem3} Let $\{e_i\}$ be a Rosenthal basis for a Banach space
$X$.
Then $\{e_i\}$ is equivalent to $c_0$ or $l_1$ or $X$ is reflexive.
In the last case, all
spreading models in
$X$ are equivalent to
$\{e_i\}$.\end{lem}

 \pf As before we can assume that $\{e_i\}$ is suppression
unconditional and invariant under spreading.
Now by James's theorem, $X$ is reflexive or contains
a subspace isomorphic to $c_0$ or $l_1$. In the last case, we may assume
$c_0$ or $l_1$ is equivalent to a block subspace of $X$, so that
$\{e_i\}$ itself is equivalent to $c_0$ or $l_1$.
If now $X$ is reflexive, any spreading model is generated by a
weakly null sequence, so by a block basic sequence, so once again
by Rosenthal's property, is equivalent to $\{e_i\}$.
\pff

Let $\{x_i\}$ and $\{y_i\}$ be basic sequences and $K>0$.  $\{x_i\}$
$K$-dominates  $\{y_i\}$ (written  $\{x_i\}\geq^K\{y_i\}$), if
$$K\|\sum a_i{x}_i\|\geq\|\sum a_i{y}_i\| \;\;\; \a (a_i)\in c_{00}$$

We will need a recent result of An\-droula\-kis, Odell, Schlum\-precht,
Tom\-czak-\-Jaegermann:

\begin{prop} {\rm (Androulakis, Odell, Schlumprecht,
Tomczak-Jaegermann)} \label{prop4}  Let $\{x_i^n\}_i, \; n\in \om$ 
be a sequence of
normalized basic weakly null sequences in a Banach space
$  X$ which have spreading models $\{\tilde{x}^n_i\}_i$, then there exists
a seminormalized basic weakly null sequence $\{y_i\}$
in
$  X$ with spreading model $\{\tilde{y}_i\}$ such that
$$2^n\|\sum a_i\tilde{y}_i\|\geq\|\sum a_i\tilde{x}_i^n\| \;\;\;\;\a n\; \a
(a_i)\in c_{00}$$
\end{prop}

\begin{cor}\label{cor5} Suppose $X$ is a Banach space with
a basis $\{e_i\}$ with Rosenthal's property. Then there exists $K>0$
such that
$\{e_i\}$ $K$-dominates any identically distributed normalized block
basis of $\{e_i\}$.\end{cor}

\pf We may assume that $\{e_i\}$ is
a normalized suppression
unconditional basic sequence, invariant under
 spreading. Assume that for any $n<\om$ there is
a normalized basic weakly null sequence $\{x_i^n\}$ in $  X$ with
 spreading model $\{\tilde{x}_i^n\}$ such that  $\{e_i\}$ does not
$4^n$-dominate  $\{\tilde{x}_i^n\}$.  Take $\{y_i\}$ as in
Proposition \ref{prop4} and notice that its spreading model  must be equivalent to
$\{e_i\}$ by Lemma \ref{lem3}. So there is some $K$ such
that $\{e_i\}\geq^K\{\tilde{y}_i\}\geq^{2^n}\{\tilde{x}_i^n\}$; taking $n$
large enough you get a contradiction. So there is
some $K$ such that $\{e_i\}$ $K$-dominates any spreading model generated by
a normalized basic
 weakly null sequence. Now any identically distributed normalized block
basis is invariant under spreading, so is its own spreading model,
giving the result. \pff

\begin{cor}\label{cor6} Suppose that $\{e_i\}$ is a normalized suppression
unconditional basic sequence, invariant under
 spreading, with Rosenthal's property. If $1$ is in 
$\{e_i\}$'s Krivine set, then $\{e_i\}\sim l_1$.\end{cor}

\pf For any $k$ take some norm one block $x(k)$ on $\{e_i\}$, such that taking
$k$ successive
 copies $x_1(k)<x_2(k)<\ldots<x_k(k)$ of it you get a sequence
$2$-equivalent to the unit vector basis in $l_1^k$.
An infinite sequence of successive copies of this vector
$x_0(k)<x_1(k)<\ldots$  must be
$K$-dominated by $\{e_i\}$, where $K$ is the constant given
by Corollary \ref{cor5}. So for any $k$,
$\{e_i\}_{i=0}^{k-1}$ $2K$-dominates $l_1^k$, but must itself, by the
triangle inequality, be $1$-dominated by $l_1^k$. Hence
$\{e_i\}$ is equivalent to $l_1$. \pff

The following proposition states that the answer to Rosenthal question
is positive if we also assume Rosenthal property in the dual.

\begin{prop}\label{prop7} Let $X$ be a  Banach space with a Rosenthal
basis, and such that  $X^*$ has a Rosenthal basis.
Then
 $X$ is isomorphic to $c_0$ or $l_p, p>1$ (and any Rosenthal basis of $X$
is equivalent to the unit vector basis
of $c_0$ or $l_p$).\end{prop}

\pf By Lemma \ref{lem3} we may assume that $X$ is reflexive.
Let $\{e_i\}$
be a Rosenthal basis of
$X$. By
renorming we may assume that the basis is suppression unconditional and
invariant under spreading. The biorthogonal basis
$\{e_i^*\}$ satisfies these
properties as well; in particular it is its own spreading model, so by
Lemma \ref{lem3}, it is equivalent to any Rosenthal basis of $X^*$; so it has
Rosenthal's property. By
Corollary \ref{cor5}, there exists
$K>0$  such that any normalized identically distributed block basis in
$X$ (resp. $X^*$) is
$K$-dominated by $\{e_n\}$ (resp. $\{e_n^*\}$).
Given a normalized identically distributed block basis $\{x_n\}$ in $X$, 
denote by $\{x_n^*\}$
a normalized identically distributed block basis in $X^*$ such
that each $x_n^*$ satisfies
$x_n^*(x_n)=1$ and has support no larger than the support of $x_n$ (this is
possible
by $1$-unconditionality and $1$-subsymmetry): $\{x_n^*\}$ is $K$-dominated
by $\{e_n^*\}$. It follows that
$\{x_n\}$ $1/K$-dominates $\{e_n\}$. Indeed, for $(a_i) \in c_{00}$,
$$K\norm{\sum a_i x_i}  \geq K\sup_{(b_i) \in c_{00}}
\frac{(\sum b_i x_i^*)(\sum a_i x_i)}{\norm{\sum b_i x_i^*}}
\geq
\sup_{(b_i) \in c_{00}} \frac{\sum b_i a_i}{\norm{\sum b_i e_i^*}}
= \norm{\sum a_i e_i}.$$
Hence any identically distributed
normalized block $\{x_i\}$ of $\{e_i\}$ is $K^2$-equi\-va\-lent to
$\{e_i\}$. By Proposition \ref{prop1}, $\{e_i\}$ must be equivalent to the
unit basis of
$l_p$ for some $p>1$. \pff

We now prove that if the selection of the subsequence in
Rosenthal's property is
continuous then the answer to the problem is also
positive; in fact we get more, it is enough to find a continuous
selection of subsequences dominating the basis. We let
$bb(\{e_i\})$ (or $bb(X)$) be the set of normalized block bases of
$\{e_i\}$, denote by
$bb_D(X)$ the same set equipped
with the product of
the discrete topology on $X$, by
$bb_E(X)$ the same set equipped with the "Ellentuck-Gowers" topology: 
basic open sets are of the form $[a,A]$ with $a<A$ for
$a=(a_1,\ldots,a_n)$ a finite normalized
block sequence and $A$ an infinite normalized block sequence,
where $$[a,A]=\{a^\frown x, x \in bb(A)\}.$$ Here $a^\frown x
\in bb(X)$ denotes the concatenation of $a$ and $x$, and
$bb(A)$ denotes the set of normalized block bases of
$A$. Proposition \ref{prop9} uses the weakest notion of continuity
combining the two topologies. We first prove a Lemma.

\begin{lem}\label{lem8} Assume $X$ is a Banach space with a Rosenthal basis
$\{e_i\}$ not equivalent to $c_0$ or $l_p$,
and let $\phi:bb(X) \rightarrow bb(X)$  map any $x \in bb(X)$ to
a subsequence of $x$. Then for any
$bb_E(X)$-open set
$[a,A]$ in
$X$ and all $n>0$, there exists a normalized block basis $x$ in $bb(A)$
such that
$\phi(a^\frown x)$  does not
$n$-dominate
$\{e_i\}$.\end{lem}

\pf Otherwise passing
to a further block, we may assume that $A=\{A_i\}$ is $C$-equivalent to
$\{e_i\}$ for some $C$; by Corollary\ref{cor5} there exists $K$ such that
$\{A_i\}$
$K$-dominates any of its identically distributed
blocks; furthermore by the assumption any block $x$ of $A$ is such
that $\phi(a^\frown x)$ $n$-dominates $\{e_i\}$,
so by $1$-subsymmetry of $\{e_i\}$, some subsequence of $x$
 $n$-dominates
$\{e_i\}$, thus $nC$-dominates $\{A_i\}$;  so
by Proposition \ref{prop1},
$A$ would be equivalent to $c_0$  or $l_p$. \pff

\begin{prop}\label{prop9} Assume $X$ is a Banach space with a
Rosenthal basis
$\{e_i\}$, and that $\phi:bb(X) \rightarrow bb(X)$ is a $bb_E(X) - bb_D(X)$
continuous map
such that for any normalized basic sequence
$x=(x_i)$ in $bb(X)$, the sequence $\phi(x)$ is a subsequence of $x$
which dominates $\{e_i\}$. Then $\{e_i\}$ is equivalent to the unit
vector basis of $c_0$ or $l_p$.\end{prop}

\pf Otherwise we build a block sequence $z=\{z_i\}$ with $\phi(z)$ not dominating
$\{e_n\}$ by induction, using Lemma \ref{lem8}.
Let $X^1=\{x^1_i\}$ be a block such that  $\phi(X^1)$ does not
$2$-dominate $\{e_i\}$. There exists an integer $N_1$ such that
$\{\phi(X^1)_i\}_{1 \leq i \leq N_1}$ does not $2$-dominate
$\{e_i\}_{1 \leq i \leq N_1}$. By continuity of $\phi$, there exists
$n_1$ in $\N$ and $A_1$ in $bb(X)$, with
 $x^1=(x^1_1,\ldots,x^1_{n_1})<A_1$, such
that if
$y=\{y_i\}$ is in
$bb(A_1)$  then
$(\phi({x^1}^\frown y))_j=(\phi(X^1))_j$ for all $1
\leq j \leq N_1$. We let $z_i=x_i^1$ for $1 \leq i \leq n_1$.

Now let $X^2=\{x^2_i\}$ be a block in $A_1$ such that
$\phi({x^1}^\frown X^2)$ does not
$4$-dominate
$\{e_i\}$. There exists an integer $N_2>N_1$ such that
the sequence $\{\phi({x^1}^\frown X^2)_i\}, 1 \leq i \leq N_2$
does not $4$ dominate
$\{e_i\}_{1 \leq i \leq N_2}$. By continuity of $\phi$,
 there exists
$n_2$ and $A_2 \in bb(A_1)$,with
 $x^2=(x^2_1,\ldots,x^2_{n_2})<A_2$ such that
if $y \in bb(A_2)$  then
$(\phi({x^1}^\frown {x^2}^\frown y))_j=
(\phi({x^1}^\frown X^2))_j$
for all $1 \leq j \leq N_2$.
We let $z_{n_1+i}=x_i^2$ for $1 \leq i \leq n_2$.
Repeating this procedure, we obtain by induction
a normalized block sequence $z=\{z_i\}$,  an
increasing sequence of
integers $\{N_i\}$, a sequence of integers $\{n_i\}$,
finite blocks $x^{i}$,  such that for
all $k$,
$$\{(\phi(z))_i\}_{1
\leq i \leq
N_k}=\{\phi({x^1}^\frown {x^2}^\frown
\ldots^\frown x^k)_i\}_{1
\leq i \leq N_k},$$
so that
$(\phi(z)_i)_{1 \leq i \leq N_k}$ does not
$2^{k}$ dominate
$\{e_i\}_{1 \leq i
\leq N_k}$, and so
$\phi(z)$ does not $2^{k}$ dominate $\{e_i\}$.
As $k$ is arbitrary
this contradicts the definition of $\phi$. \pff

Let us  remark, that once again, there is an opposition between
Rosenthal's property and block equivalence minimality. Indeed, for any
block equivalence minimal basis, Gowers' theorem implies the existence of
a winning strategy to produce block sequences ($C$-) equivalent to 
$\{e_i\}$; in Gowers' game defined by Bagaria and Lopez-Abad, 
which is actually equivalent to the original game defined by 
Gowers (\cite{BL}), Player 1 plays block vectors
and Player 2 sometimes chooses a vector in the finite dimensional space defined
 by the blocks played by Player 1.
The winning strategy then defines a continuous map from block
sequences to further block sequences ($C$-) equivalent to
$\{e_i\}$. Notice also that Gowers-Maurey constructions (\cite{gm}) yield 
winning strategies in the
previous sense: technically, $l_1^n$-averages used to build interesting vectors
in their space may at each step of the construction
be chosen in an arbitrary block-subspace. Roughly speaking, this
means that if one
tried to adapt their ideas to build a non-trivial Rosenthal basis,
not only one would have to find a way to pass 
from selecting further (finite) blocks to  selecting (infinite)
subsequences, but also one would probably have to add new methods to suppress 
the continuity of the selection map.

\

Finally, we investigate the relation between
Rosenthal's question and a problem of S. Argyros.

\paragraph{Problem} (S. Argyros) {\em Let $X$ be a Banach space such that all
spreading models in $X$ are equivalent. Must these
spreading models be equivalent to the unit vector basis or
$l_p$
for  some $p \geq 1$?}

\

For example, spaces $l_p$ have unique spreading model up to
equivalence. Indeed, in the reflexive case, all spreading models
are generated by weakly null sequences; and in $l_1$, any
spreading model is generated by a
$l_1$-sequence, or by Rosenthal's theorem, by a weakly Cauchy sequence.
In the second case, the difference sequence is weakly null, so generates
$l_1$, and it follows that the spreading model is equivalent to $l_1$.
But this does not generalize to the case of $c_0$, since the unit
basis of $c_0$ and the summing basis generate non-equivalent
spreading models; however all spreading models generated by weakly
null sequences are clearly equivalent to $c_0$.

\

Lemma \ref{lem3} shows that a  positive answer to the problem of Argyros implies
a positive answer to the problem of Rosenthal. Actually
Androulakis, Odell,
Schlumprecht and Tomczak-Jaegermann proved that the
answer to Argyros' Problem is positive under the additional assumption of
uniformity or that $1$ is in the Krivine set
of some basic sequence. Our methods are inspired from their results.
A natural generalization of Argyros' question is mentioned in their
article: if a Banach space contains only countably many
spreading models up to equivalence, must one of them be equivalent
to $c_0$ or $l_p$? In the other direction, the
following remark about Banach spaces with more than  countably many
spreading models is a straightforward consequence of a
well-known result of Silver.

\begin{prop}\label{prop10} Let $X$ be a separable
Banach space. Then either
$X$ contains continuum many non-equivalent spreading models, or $X$
contains at most countably many non-equivalent spreading models.
When
$X^*$ is separable, the
same dichotomy holds for spreading models generated by weakly null
basic sequences; when $X$ has a Schauder basis, it holds for
spreading models generated by block basic sequences.\end{prop}

\pf In the following, $\sim^C$ denotes the usual
$C$-equivalence between basic sequences.
We consider the set $\cal S$ of semi-normalized basic sequences
generating spreading models, which can be
described as
the set of semi-normalized basic sequences $\{x_i\}$ such that:
$\forall k \in \N,\forall \eps>0, \exists N: \forall N<l_0<\ldots<l_k,
\forall N<l^{\prime}_0<\ldots<l^{\prime}_k,
(x_{l_i})_{i=0}^k
\sim^{1+\eps} (x_{l'_i})_{i=0}^k.$
This  set is clearly a Borel subset of the Polish
space $X^{\omega}$.
Now consider the equivalence relation $\simeq$ on $\cal S$ meaning that the two
sequences generate
spreading models which are equivalent in the usual $\sim$ sense. That is
$(y_n) \simeq (z_n)$ iff
$$\exists C>0, \forall k \in \N,\exists N: \forall N<l_0<\ldots<l_k,
(y_{l_i})_{i=0}^k \sim^C (z_{l_i})_{i=0}^k.$$ This equivalence
relation is Borel as well. Now by a Theorem of Silver (Th 35.20 in
\cite{kec}), a Borel (even coanalytic) equivalence relation on a
Borel subset of a Polish space has either only countably many
classes or there exists a Cantor set of mutually non-equivalent
elements. As two spreading models are $\sim$-equivalent if and
only if any two semi-normalized basic sequences which generate
them are
$\simeq$-equivalent, the result follows.
When $X^*$ is separable, the same proof holds for the set of
weakly null semi-normalized basic sequences generating spreading
models, which is also Borel in $X^{\omega}$; or when $X$ has a basis, for
the set of block basic sequences generating spreading models.\pff

\paragraph{Remark} It is also a consequence of the Theorem of Silver
that a Schauder basis of a Banach space $X$ has continuum many
non-equivalent subsymmetric block basic sequences or only
countably many classes of equivalence of them. Indeed, the set of
of normalized block bases equipped with the product topology on
$X$ is Polish, and the set of subsymmetric normalized block basic
sequences is $F_{\sigma}$ in it.

\

A result  analogous to Proposition \ref{prop9} holds also for spreading models.
It turns out that the continuity of a map, which picks subsequences generating 
spreading models in a strong sense described below, 
is strong enough to imply that
there is actually a copy of $\ell_p$ or $c_0$ in the space.

First some terminology: given a sequence
$\eps=\{ \eps_i\}$, $\eps_i\searrow 0$, and a basic sequence 
$\{x_i\}$ in a Banach space we say that $\{x_i\}$ 
$\eps$-generates a spreading model $\{\tilde{x}_i\}$, if 
for any $k<n_1<n_2<\dots <n_k$ we have 
$(x_{n_1},\dots ,x_{n_k})\sim^{1+\eps_k}(\tilde{x}_1,\dots ,\tilde{x}_k)$.
Obviously every basic sequence for any sequence $\eps$ of non-zero 
scalars  has a subsequence 
$\eps$-generating a spreading model. We will use the notation introduced 
before Lemma \ref{lem8}.

We say that a Banach space 
$X$ contains almost isometric copies of $c_0$ (resp. $\ell_p$), if 
for any $\delta >0$, $X$ has a subspace $(1+\delta)$-isomorphic to 
$c_0$ (resp. $\ell_p$). 

\begin{prop}\label{prop11} Let $X$ be a Banach space with a basis $\{e_i\}$. 
Fix a sequence $\eps=\{\eps_i\}$, $\eps_i\searrow 0$.
 Assume there is a continuous map $\phi:bb_D(X) \rightarrow bb_D(X)$ 
such that for any normalized basic sequence $x=\{ x_i\}$ 
in $bb(X)$, the sequence $\phi(x)$ is a subsequence of $x$ 
which $\eps$-generates a spreading model. Then $X$ contains  
almost isometric copies of $c_0$ or $l_p$ for some $1\leq p<\infty$.\end{prop}
 
\pf We recall the notion of asymptotic spaces as presented in 
\cite{ode}. Let $X$ be a Banach space with a basis $\{ e_i\}$. 
A tail subspace means here a block subspace of $X$ of a finite 
codimension. 

We say that a normalized 
basic sequence $\{ a_i\}_{i=1}^{n}$ is asymptotic in $X$, 
if 
$$\forall \delta>0\;\;\forall k_1\;\;\exists x_1\in\langle e_i\rangle_{i>k_1}\;\;
\forall k_2\;\;\exists x_2\in\langle e_i\rangle_{i>k_2}\;\dots\;
\forall k_n\;\;\exists x_n\in\langle e_i\rangle_{i>k_n}$$
so that $\{x_i\}_{i=1}^n$ is a normalized block sequence 
$(1+\delta )$-equivalent to $\{ a_i\}_{i=1}^n$.
In other words, if we consider the asymptotic game, in which player I 
picks tail 
subspaces and player II picks block vectors from the 
subspaces chosen by player 
I, then a normalized basic sequence $\{a_i\}_{i=1}^n$ is asymptotic iff 
player II for any $\delta$  has a winning strategy in choosing a 
normalized block sequence of vectors $(1+\delta)$-equivalent to 
$\{a_i\}_{i=1}^n$. 

Since a block sequence has a subsequence generating an unconditional 
spreading model, by Krivine's theorem, there is some $1\leq p\leq \infty$ 
such that $\ell_p^n$ or $c_0^n$ (in case $p=\infty$) is asymptotic
in $X$ for any $n\in\N$ (i.e. the unit basic vectors in 
$\ell_p^n$ or $c_0^n$ form asymptotic sequences). 

Let $X$ satisfy the assumption of the proposition. We will use in the proof
only asymptotic sequences of length 2. Pick $1\leq p\leq \infty$
such that $\ell_p^2$ or $c_0^2$ (in case $p=\infty$) is asymptotic
in $X$. 

We will show that any asymptotic pair  
$(a_1,a_2)$ of $X$ is 1-equivalent to 
the unit basic vectors of suitable $\ell_p^2$ or $c_0^2$.

Fix $\delta >0$ and pick any asymptotic pair 
$(a_1,a_2)$. Pick $n\in\N$, $n>1$ such 
that $(1+\eps_{n-1})^3<1+\delta$.
Consider the asymptotic game for $(1+\eps_{n-1})$ and $(a_1,a_2)$. 
Let $x=(x_1,x_2,\dots )$ be a block sequence consisting of 
vectors picked by player II 
in the first move in some game (you can produce such a sequence by 
letting player I choose in the first move tail subspaces of arbitrary 
large codimension). Let $\phi (x)=(x_{j_1},x_{j_2},\dots)$. 
By the continuity of $\phi$ there is some 
$J>j_n$ such that 
$\phi ([(x_1,\dots , x_{J}),X])\subset [(x_{j_1},
\dots , x_{j_n}), X]$.
 
Now consider sequence $y=(x_1,\dots x_{J}, y_1,y_2,\dots)$, where 
$ (y_1, y_2,\dots )$, with $x_J<y_1$, is a block sequence of vectors 
chosen by player II in 
the second move in some game for $(1+\eps_{n-1})$ and $(a_1,a_2)$, 
in which player II picked in the first move the vector
$x_{j_n}$ (again you produce such sequence by letting player I choose in the
second move 
tail subspaces of arbitrary large codimension). Let 
$$\phi (y)=(x_{j_1},\dots, x_{j_n},\dots x_{j'}, y_{k_1}, 
y_{k_2}, \dots)$$ 
Again by the continuity of $\phi$ there is some $K>k_1$ such that 
$$\phi ([(x_1,\dots ,x_{J}, y_1,\dots ,y_{K}),X])
\subset [(x_{j_1},\dots, x_{j_n},\dots, x_{j'}, y_{k_1}),X].$$ 

Now consider the asymptotic game for $(1+\eps_{n-1})$ and asymptotic 
$\ell_p^2$ or $c_0^2$. Repeating the previous 
procedure for such $\ell_p$ or $c_0$ we extend the finite sequence 
$(x_1,\dots ,x_{J}, y_1,\dots ,y_{K})$ by suitable block sequences 
$(1+\eps_{n-1})$-realizing 
$\ell_p^2$ or $c_0^2$ in $X$. In this way we obtain finite block sequences
$$b=(x_1,\dots, x_{J},y_1,\dots y_{K},v_1,\dots v_{L},
z_1,\dots, z_{M})$$ and
$$c=(x_{j_1},\dots,x_{j_n},\dots ,x_{j'},y_{k_1},\dots ,y_{k'},
v_{l_1},\dots ,v_{l'},z_{m_1})$$
with 
$(x_{j_n},y_{k_1})\sim^{1+\eps_{n-1}}(a_1,a_2)$, 
$(v_{l_1},z_{m_1})\sim^{1+\eps_{n-1}}\ell_p$ or $c_0$,
and $\phi ([b,X])\subset [c,X]$. 

By definition of $\phi$, 
$(x_{j_n},y_{k_1})\sim^{1+\eps_{n-1}}(v_{l_1},z_{m_1})$,
Hence, by the choice of $n$, we have 
$(a_1,a_2)\sim^{1+\delta}\ell_p^2$ or $c_0^2$. 
Since $\delta$ was arbitrary small, $(a_1,a_2)$ is 
1-equivalent to the unit vector basis of $\ell_p^2$ or $c_0^2$. 

Now by a standard procedure we produce a 
block subspace $X_0$ of $X$ such that for any 
$\delta >0$ there is a $N_{\delta}$ such that for any $N_{\delta}<x_1<x_2$, 
$x_1,x_2\in X_0$ we have $(x_1,x_2)\sim^{(1+\delta)}\ell_p^2$ or $c_0^2$, and 
then produce an almost isometric copy of $\ell_p$ or $c_0$ in the space 
(cf. e.g. \cite{mit}, \cite{mmt}), which finishes the proof. \pff

\paragraph{Remark} Notice that by this proposition 
in $\ell_p$, for $1<p<\infty$, endowed 
with a distorting norm one cannot pick sequences producing 
spreading models (in 
the sense defined above) in a continuous way, however any block sequence is
subsymmetric.

\

\paragraph{Acknowledgements} We wish to thank G. Androu\-lakis and
T.Schlum\-precht for useful
informa\-tion about spreading models and comments about this paper.

\

\small\em

Valentin Ferenczi, Christian Rosendal 

Equipe d'Analyse, Boite 186, Universit\'e Paris 6

4, place Jussieu, 75252 Paris Cedex 05, France

e-mail: ferenczi@ccr.jussieu.fr, rosendal@ccr.jussieu.fr

\

Anna Pelczar 

Jagiellonian University

Reymonta 4, 30-059 Krak\'ow, Poland

e-mail: apelczar@im.uj.edu.pl

\end{document}